\documentclass[12pt,twoside,a4paper]{amsart}
\usepackage{amssymb,amsfonts,amsmath}
\usepackage[all,knot,arc]{xy}
\usepackage[dvips]{graphicx}
\usepackage{psfrag}

\def\endpf{\relax\ifmmode\expandafter\endproofmath\else
  \unskip\nobreak\hfil\penalty50\hskip.75em\hbox{}\nobreak\hfil\bull
  {\parfillskip=0pt \finalhyphendemerits=0 \bigbreak}\fi}
\def\bull{\vbox{\hrule\hbox{\vrule\kern3pt\vbox{\kern6pt}\kern3pt\vrule}\hrule}}

\newtheorem{defn}{Definition}[section]
\newtheorem{lemma}[defn]{Lemma}

\newtheorem{proposition}[defn]{Proposition}

\newtheorem{maintheorem}{Theorem}

\newtheorem{maincor}[maintheorem]{Corollary}

\newcommand{\zz}{{\mathbb Z}}
\newcommand{\rr}{{\mathbb R}}

\newcommand{\qq}{{\mathbb Q}}

\newcommand{\CC}{{\mathcal C}}
\newcommand{\cp}{{\mathbb C}{\mathrm P}^2}

\newcommand{\poincare}{Poincar\'{e}}

\newcommand{\ozsvath}{Ozsv\'{a}th}
\newcommand{\szabo}{Szab\'{o}}
\newcommand{\spin}{\ifmmode{\rm Spin}\else{${\rm spin}$\ }\fi}
\newcommand{\spinc}{\ifmmode{{\rm Spin}^c}\else{${\rm spin}^c$\ }\fi}

\newcommand{\tors}{{\it Tors}}

\newcommand{\ds}{\displaystyle}
\def\ceil#1{\left\lceil#1\right\rceil}

\DeclareMathOperator{\lk}{lk}

\newcommand{\Tpq}{T_{p,q}}

\def\srput(#1)#2{\rput(#1){\tiny #2}}

\newenvironment{narrow}[2]{%
 \begin{list}{}{%
  \setlength{\topsep}{0pt}%
  \setlength{\leftmargin}{#1}%
  \setlength{\rightmargin}{#2}%
  \setlength{\listparindent}{\parindent}%
  \setlength{\itemindent}{\parindent}%
  \setlength{\parsep}{\parskip}%
 }%
\item[]}{\end{list}}

\newif\ifpic
\picfalse    
\pictrue   

\DeclareMathOperator\Hom{Hom}

\begin{document}

\title{Dehn surgeries and negative-definite 
four-manifolds}
\author[Brendan Owens]{Brendan Owens}
\address{School of Mathematics and Statistics \newline\indent 
University of Glasgow \newline\indent 
Glasgow, G12 8QW, United Kingdom}
\email{brendan.owens@glasgow.ac.uk}
\author[Sa\v so Strle]{Sa\v so Strle}
\address{Faculty of Mathematics and Physics \newline\indent 
University of Ljubljana \newline\indent Jadranska 21 \newline\indent 
1000 Ljubljana, Slovenia }
\email{saso.strle@fmf.uni-lj.si}
\date{\today}
\thanks{B. Owens was supported in part by NSF grant DMS-0604876 and by the EPSRC.\\
S.Strle was  supported in part by the ARRS of the
Republic of Slovenia research program No. P1-0292-0101.  We also acknowledge support from the ESF through the ITGP programme}
\dedicatory{Dedicated to Jos{\'e} Maria Montesinos on the occasion of his 65th birthday.}
\begin{abstract}
Given a knot $K$ in the three-sphere, we address the question: which  Dehn surgeries on $K$ bound negative-definite four-manifolds?
We show that the answer depends on a  number $m(K)$, which is a smooth concordance invariant.
We study the properties of this invariant, and compute it for torus knots.
\end{abstract}

\maketitle

\pagestyle{myheadings}
\markboth{BRENDAN OWENS AND SA\v{S}O STRLE}{DEHN SURGERIES AND NEGATIVE-DEFINITE FOUR-MANIFOLDS}


\section{Introduction}
\label{sec:intro}

The intersection pairing of a smooth compact four-manifold, possibly with boundary, is
an integral symmetric bilinear form $Q_X$ on $H_2(X;\zz)/\tors$; 
it is nondegenerate if the boundary of $X$ is a rational homology sphere.
Given a rational homology three-sphere $Y$ there are various gauge-theoretic constraints on which bilinear forms
may be intersection pairings of manifolds bounded by $Y$.
For example, Donaldson's celebrated Theorem A \cite{d} tells us that the only negative-definite
pairings bounded by the three-sphere are the standard diagonal forms.  Another well-known example is the \poincare\
homology sphere $P$, oriented as the boundary of the positive-definite $E8$ plumbing; this does not bound any negative-definite four-manifold.  An alternative description of $P$ is $+1$ surgery on the right-handed trefoil knot.  It is also well-known that $+5$ surgery on the same trefoil knot gives the lens space $L(5,1)$ which is the boundary of a negative-definite disk bundle over $S^2$.  A natural question arises, for the trefoil and more generally for any knot $K$ in $S^3$: for which  rational numbers $r$ does the Dehn surgery $S^3_r(K)$ bound a smooth negative-definite four-manifold?  This question is related to the computation of unknotting numbers, and also to the classification of tight contact structures.

An easy argument, for example using Lemma \ref{lem:sigAn}, shows that $S^3_r(K)$  bounds negative-definite whenever $r$ is negative.  Any knot can be converted to the unknot by a finite number of crossing changes.  This enables us to show that large positive  surgeries on a knot $K$ always bound negative-definite four-manifolds, and leads to the following invariant:
\begin{equation}\label{eq:mdef}
m(K)=\ds\inf\{r\in\qq_{>0}\,\,|\,S^3_r(K)\ \mbox{bounds a negative-definite 4-manifold}\}.
\end{equation}
Some properties of $m$ are described in the following theorem which we prove in Section \ref{sec:basics}; we conjecture that the subadditivity property of the last part holds for $m(K)$ as well.

\begin{maintheorem}\label{thm:properties}
(a) Let $K \subset S^3$ be a knot. If $K$ can be unknotted by changing $p$ positive  and $n$ negative crossings, then $m(K)\le 4p$.  In particular, if $K$ admits a diagram without positive crossings then $m(K)=0$.\\
(b) For all rational numbers $r>m(K)$, the $r$-surgery $S^3_r(K)$ on $K$ bounds a negative-definite manifold $X_r$ with $H^2(X_r)\to H^2(S^3_r(K))$ surjective.\\
(c) $m(K)$ is a concordance invariant of $K$, hence it defines a function $m:\CC \to \rr_{\ge0}$, where $\CC$ denotes the smooth concordance group of classical knots.\\
(d) The integer valued invariant $\ceil m$ is subadditive with respect to connected sum, i.e.~
$$\ceil {m(K\#C)}\le\ceil {m(K)}+\ceil{m(C)}$$
for any knots $K$ and $C$ in $S^3$.
\end{maintheorem}

In Section \ref{sec:torus} we compute $m$ for torus knots.  The question of which nonzero Dehn surgeries on torus knots bound negative-definite manifolds was considered previously in \cite{short} and \cite{greene}; we give a complete answer here.
Let $p$, $q$ be coprime integers with $p>q>0$.  It is well known that the $(pq-1)$-surgery on the torus knot $\Tpq$ bounds a negative-definite manifold (since this surgery is a lens space \cite{moser}). 
Let $n$ be the number of steps in the standard Euclidean algorithm for $p$ and $q$.
Denote by $q^*=q^{-1}\pmod p$ the solution to the congruence $qa\equiv1\pmod p$ with $0<a<p$, and similarly let $p^*=p^{-1}\pmod q$. 
\begin{maintheorem}\label{thm:Tpq}
Let $\Tpq$ denote the positive $(p,q)$-torus knot.  Then 
$$m(\Tpq)=
\left\{
\begin{array}{ll}
\ds pq-\frac{q}{p^*} & \quad\mbox{if $n$ is even,}\\
\ds pq-\frac{p}{q^*} & \quad\mbox{if $n$ is odd.}\\
\end{array}
\right.$$
The manifold given by $m(\Tpq)$ surgery on $\Tpq$ bounds a negative-definite four-manifold. Moreover, for any negative-definite four-manifold this surgery bounds not all $\spinc$ structures on the surgery manifold extend over the four-manifold.

For any negative torus knot, $m(T_{p,-q})=0$.
\end{maintheorem}

The special case of $q=n=2$ follows from work of Lisca-Stipsicz \cite{ls}; they used essentially the same obstruction to get a lower bound on $m$ but a different construction for the negative-definite manifold bounded by $S^3_{m(\Tpq)}(\Tpq)$.  Theorem \ref{thm:Tpq} yields the following generalisation of \cite[Theorem 4.2]{ls}.

\begin{maincor}\label{cor:Tpq}
For each rational number $r$ in the interval $[pq-p-q, m(\Tpq))$, the 3-manifold given by $r$ surgery on $\Tpq$ does not admit fillable contact structures.
\end{maincor}


\section{Continued fractions and surgery cobordisms}
\label{sec:contfrac}
In this section we establish some notation and basic facts about continued fractions and surgery cobordisms.
Given a sequence of numbers $c_1, c_2, \dots, c_n$ in the extended real line $\rr\cup\{\infty\}$ (typically these will be integers) one obtains two numbers
$$[c_1, c_2, \dots, c_n]^+:=c_1+\frac{1}{c_2+\raisebox{-3mm}{$\ddots$
\raisebox{-2mm}{${+\frac{1}{\displaystyle{c_n}}}$}}}\, ,$$
and
$$[c_1, c_2, \dots, c_n]^-:=c_1-\frac{1}{c_2-\raisebox{-3mm}{$\ddots$
\raisebox{-2mm}{${-\frac{1}{\displaystyle{c_n}}}$}}}\, ;$$
we refer to these as positive and negative continued fraction expansions, respectively.  
For a given positive rational number the standard Euclidean algorithm yields a unique expansion
$$\frac pq=[c_1, c_2, \dots, c_n]^+$$
with integer coefficients satisfying $c_1\ge0$, $c_i>0$ for $1<i<n$ and $c_n>1$.
Similarly there is a unique expansion
$$\frac pq=[a_1, a_2, \dots, a_m]^-,$$
with positive integer coefficients and $a_i>1$ for all $i>1$.  (This is often referred to as the Hirzebruch-Jung continued fraction of $\frac pq$.)

The following is immediate from the Euclidean algorithm for $p/q$.

\begin{lemma}
\label{lem:cf0}
Let $p>q>1$ be coprime integers with 
$$\frac pq=[c_1, c_2, \dots, c_n]^+.$$  
Writing $p=c_1q+r$ we have
$$\frac qr=[c_2, \dots, c_n]^+.$$
\end{lemma}

\begin{lemma}
\label{lem:cf1}
The positive and negative continued fraction expansions are related by the formula
$$[c_1,\ldots,c_n]^+=[c_1+1,\underbrace{2,\ldots,2}_{c_2-1},c_3+2,\underbrace{2,\ldots,2}_{c_4-1},c_5+2,\ldots ]^-,$$
where the continued fraction ends with $c_n+1$ if $n$ is odd and with $c_{n}-1$ $2$'s if $n$ is even.
\end{lemma}

\proof
This may be deduced easily from the equations
$$[c,x]^+=\left[c+1,\frac{x}{x-1}\right]^-$$
and
\begin{equation}\label{eqn:cf}
y=[c,z]^+\iff \ds\frac y{y-1}=[\underbrace{2,\ldots,2}_{c-1},1+z]^-.
\end{equation}
(Or see \cite[Proposition 2.3]{ppp}.)
\endproof

\begin{lemma}
\label{lem:cf2}
Let $p>q>0$ be integers with $\dfrac pq=[c_1, c_2, \dots, c_n]^+$.  Then
\begin{equation}
\frac p{p-q}=[\underbrace{2,\ldots,2}_{c_1-1},c_2+2,\underbrace{2,\ldots,2}_{c_3-1},c_4+2,\ldots ]^-,
\end{equation}
where the continued fraction ends with $c_n+1$ if $n$ is even and with $c_{n}-1$ $2$'s if $n$ is odd.
\end{lemma}
\proof
This follows easily from Lemma \ref{lem:cf1} and \eqref{eqn:cf}.  (Or see \cite[Lemma 7.2]{neumann} or \cite[Proposition 2.7]{ppp}.)
\endproof

Given coprime natural numbers $p$ and $q$, define $q^*$ to be the multiplicative inverse of $q$ modulo $p$, i.e. $qq^*\equiv 1\pmod p$ and $1\le q^*<p$.

\begin{lemma}
\label{lem:cf3}
If $p>q$ are coprime positive integers and
$$\frac pq=[c_1,\ldots,c_n]^-$$
with integers $c_i\ge2$, then
$$\frac p{q^*}=[c_n,\ldots,c_1]^-.$$
\end{lemma}
\proof
This can be seen by induction on $n$. For details see for example \cite{hnk}.
\endproof

\begin{figure}[htbp] 
\begin{center}
\ifpic
\psfrag{p q}{\scriptsize$\ds\frac{p}{q}$}
\psfrag{c 1}{\scriptsize$c_1$}
\psfrag{c 2}{\scriptsize$c_2$}
\psfrag{c 3}{\scriptsize$c_3$}
\psfrag{c 4}{\scriptsize$c_{n-1}$}
\psfrag{c 5}{\scriptsize$c_n$}
\psfrag{equiv}{$\sim$}
\psfrag{dots}{$\dots$}
\includegraphics[width=10cm]{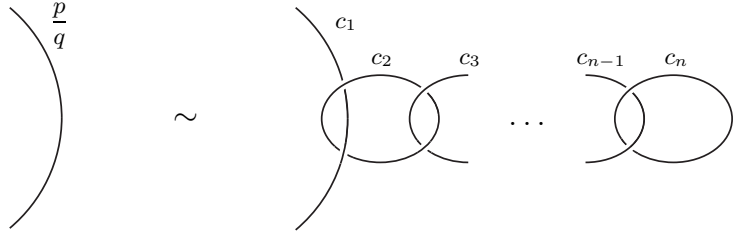}
\else \vskip 5cm \fi
\begin{narrow}{0.3in}{0.3in}
\caption{
\bf{Converting between Dehn surgery and integral surgery.}}
\label{fig:qzsurg}
\end{narrow}
\end{center}
\end{figure}

The importance of continued fractions in our context is due to their appearance when Dehn surgeries are converted to integer surgeries.  If
$$\frac pq=[c_1,\ldots,c_n]^-,$$
then the 3-manifold given by $\ds\frac pq$ surgery on a knot $K$ is equivalent to that given by a framed link consisting of $K$ with framing $c_1$ and a chain of unknots with framings $c_2,\dots,c_n$ as in Figure \ref{fig:qzsurg}. The equivalence of the two descriptions is established using the slam-dunk move (see \cite[\S 5.3]{gs} for details). Note that for surgeries on the unknot the equality of Lemma \ref{lem:cf3} combined with the integer surgery presentation corresponds to the equality $L(p,q)=L(p,q^*)$.
The next lemma is of use in computing the signature of the resulting $4$-manifold with boundary.

\begin{lemma}\label{lem:sigAn}
Let $a_1,\dots,a_n$ be integers with $a_1\ge1$, $|a_i|\ge2$ for $1<i<n$ and either $|a_n|\ge2$ or $a_n=-1$.  Let $A$ denote the symmetric $n\times n$ matrix whose nonzero entries are $A_{i,i}=a_i$, $A_{i,i\pm1}=1$. Then
$$\mathrm{signature}(A)=\#\{i \, |\,  a_i>0\} - \#\{i \, |\,  a_i<0\}.$$
\end{lemma}

\proof
Write $A_n$ for $A$. We prove the formula by induction. Let $A_{n-1}$ denote the minor given by deleting the last row and column of $A_n$. We claim that
$$\mathrm{signature}(A_n)=\mathrm{signature}(A_{n-1})+\mathrm{sign}(a_n).$$
Let $v_1,\dots,v_n$ be basis vectors for a free abelian group with a bilinear pairing given by $Q(v_i,v_j)=A_{i,j}$. Extending coefficients to $\qq$ there are constants $b_i\in\qq$ for which $\ds v'_n=v_n-\sum_{i=1}^{n-1} b_iv_i$ is $Q$ orthogonal to $v_i$ for $i<n$; in fact, $b_{n-1}=1/[a_{n-1},\ldots,a_1]^-$.  In the basis $v_1,\dots,v_{n-1},v'_n$ the form $Q$ has a block matrix with blocks $A_{n-1}$, $Q(v'_n,v'_n)$ from which it follows that the signature of $A$ is given by the sum of the signature of $A_{n-1}$ and the sign of $Q(v'_n,v'_n)$.  The conditions on $a_1,\dots,a_{n-1}$ ensure  that $-1<1/[a_{n-1},\ldots,a_1]^- \le 1$ and
$$\mathrm{sign}(Q(v'_n,v'_n))=\mathrm{sign}\left(a_n-\frac{1}{[a_{n-1},\ldots,a_1]^-}\right)=\mathrm{sign}(a_n).$$

Alternatively, note that successively blowing down $\pm 1$ entries on the diagonal we obtain a diagonally dominant matrix for which the signs of eigenvalues are given by the signs of the diagonal entries.
\endproof

We start our study of negative-definite cobordisms that determine the behaviour of $m(K)$ by showing that if some surgery on $K$ bounds a negative-definite manifold, then so does any larger surgery.

\begin{lemma}\label{lem:rscob}
Let $K$ be a knot in $S^3$ and let $r, s$ be rational numbers with $r>s>0$.  Then there exists a negative-definite two-handle cobordism from $S^3_s(K)$ to $S^3_r(K)$.
\end{lemma}
\proof
Suppose that the negative continued fractions of $r, s$ agree for the first $m$ terms, $m\ge0$.  In other words we have
\begin{eqnarray*}
s&=&[c_1,\dots,c_m,c_{m+1},\dots,c_{m+k}]^-\\
r&=&[c_1,\dots,c_m,c'_{m+1},\dots,c'_{m+k'}]^-\\
\end{eqnarray*}
with $c_1,c'_1\ge1$, $c_n,c'_n\ge2$ for all $n\ge2$, and $0<l=c'_{m+1}-c_{m+1}$.

The cobordism from $S^3_s(K)$ to $S^3_r(K)$ is then a composition of cobordisms $W_0,\dots,W_l$ which we proceed to describe:

\addtolength{\arraycolsep}{-3pt}
$$\begin{array}{cccccclc}
S^3_s(K)\stackrel{W_0}{\longrightarrow}&S^3_{[c_1,\dots,c_{m+1}]^-}(K)&\stackrel{W_1}{\longrightarrow}&S^3_{[c_1,\dots,c_{m+1}+1]^-}(K)&\stackrel{W_2}{\longrightarrow}&S^3_{[c_1,\dots,c_{m+1}+2]^-}(K)&\stackrel{W_3}{\longrightarrow}\\
&&&\cdots&\stackrel{W_{l-1}}{\longrightarrow}&S^3_{[c_1,\dots,c'_{m+1}-1]^-}(K)&\stackrel{W_l}{\longrightarrow}&S^3_r(K);\\
\end{array}$$
the negative-definiteness of each of $W_0,\dots,W_l$ follows by Lemma \ref{lem:sigAn}.

To obtain the cobordism $W_0$ note that $S^3_{[c_1,\dots,c_{m+1}]^-}(K)$ bounds the positive-definite integer surgery cobordism given by $K$ with a chain of linked unknots in the usual way, with framings given by the continued fraction coefficients.  There is an obvious positive-definite cobordism from this to the corresponding integer surgery description of $S^3_s(K)$; reversing orientation yields $W_0$.

Each $W_i$ for $0<i<l$ is the surgery cobordism given by attaching a $(-1)$-framed unknot to the positive-definite integer surgery description of $S^3_{[c_1,\dots,c_{m+1}+i-1]^-}(K)$ along the meridian of the last unknot in the chain.

Now let $r'=[c'_{m+2},\dots,c'_{m+k'}]^-$ (in other words $r'$ is given by the tail of the continued fraction of $r$) and let
$$\frac{r'}{r'-1}=[a_1,\dots,a_n]^-,$$
with $a_i\ge2$.  Then we have
\begin{eqnarray*}
r&=&[c_1,\dots,c_m,c'_{m+1},r']^-\\
&=&[c_1,\dots,c_m,c'_{m+1}-1,-a_1,-a_2,\dots,-a_n]^-,
\end{eqnarray*}
which yields the negative-definite surgery cobordism $W_l$ from $S^3_{[c_1,\dots,c_m,c'_{m+1}-1]^-}(K)$ to $S^3_r(K)$.
\endproof

The next two lemmas exhibit negative-definite cobordisms from the disjoint union of $S^3_r(K)$ and $S^3_s(C)$ to $S^3_{r+s}(K\#C)$ for certain surgery coefficients $r$ and $s$. We use them to prove subadditivity of an integer version of $m(K)$ under connected sums.  We can exhibit many other such cobordisms and in fact conjecture that such a cobordism exists for any positive rational numbers $r$ and $s$. The point of the second lemma is that we do not know if the surgery manifold $S^3_{\ceil{m(K)}}(K)$  bounds a negative-definite manifold in general.

\begin{lemma}
\label{lem:sum1}
Let $K$ and $C$ be knots in $S^3$ and let $r, s$ be positive rational numbers whose sum is an integer.  Then there exists a negative-definite cobordism from $S^3_r(K)\sqcup S^3_s(C)$ to $S^3_{r+s}(K\#C)$.
\end{lemma}

\begin{figure}[htbp] 
\begin{center}
\ifpic
\psfrag{K}{\scriptsize$K$}
\psfrag{C}{\scriptsize$C$}
\psfrag{m}{\scriptsize$m$}
\psfrag{n}{\scriptsize$n$}
\psfrag{0}{\scriptsize$0$}
\psfrag{arrow}{$\longrightarrow$}
\includegraphics[width=10cm]{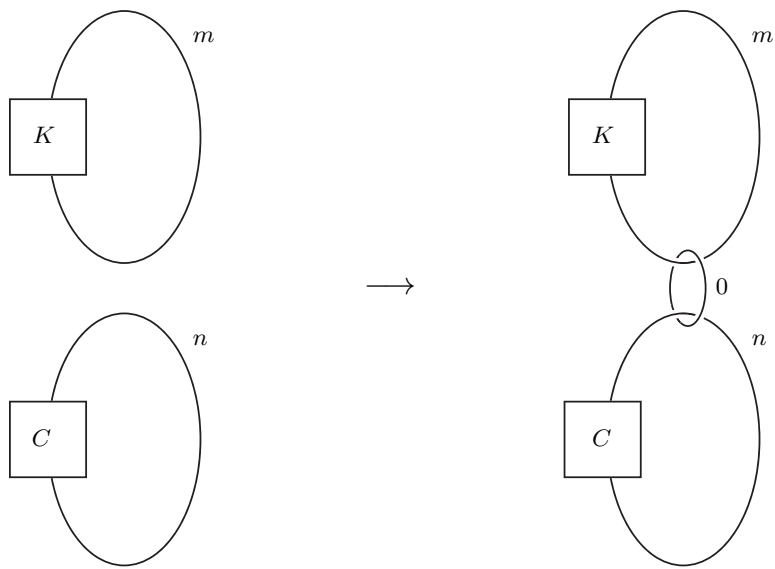}
\else \vskip 5cm \fi
\begin{narrow}{0.3in}{0.3in}
\caption{
\bf{A negative-definite cobordism from $S^3_m(K)\sqcup S^3_n(C)$ to $S^3_{m+n}(K\#C)$.}}
\label{fig:sum1}
\end{narrow}
\end{center}
\end{figure}

\proof
We begin with the simplest case which is when $r=m$ and $s=n$ are both integers.  The Kirby diagram on the right hand side of Figure \ref{fig:sum1} represents a four-manifold whose boundary is $S^3_{m+n}(K\#C)$: to see this trade the 0-framed 2-handle for a 1-handle in dotted circle notation, slide one of the remaining 2-handles over the other and cancel the 1-handle with a 2-handle.  The addition of the 0-framed 2-handle gives the required negative-definite cobordism in this case, since the four-manifold on the left of   Figure \ref{fig:sum1} is definite with $b_2^+=2$ and that on the right has $b_2^+=2$ and $b_2^-=1$.
(Similar calculations show the cobordism is also negative-definite if one or both of $m, n$ is zero.)

\begin{figure}[htbp] 
\begin{center}
\ifpic
\psfrag{K}{\scriptsize$K$}
\psfrag{C}{\scriptsize$C$}
\psfrag{m}{\scriptsize$m$}
\psfrag{n}{\scriptsize$n$}
\psfrag{0}{\scriptsize$0$}
\psfrag{1}{\scriptsize$1$}
\psfrag{2}{\scriptsize$2$}
\psfrag{3}{\scriptsize$3$}
\psfrag{arrow}{$\longrightarrow$}
\includegraphics[width=10cm]{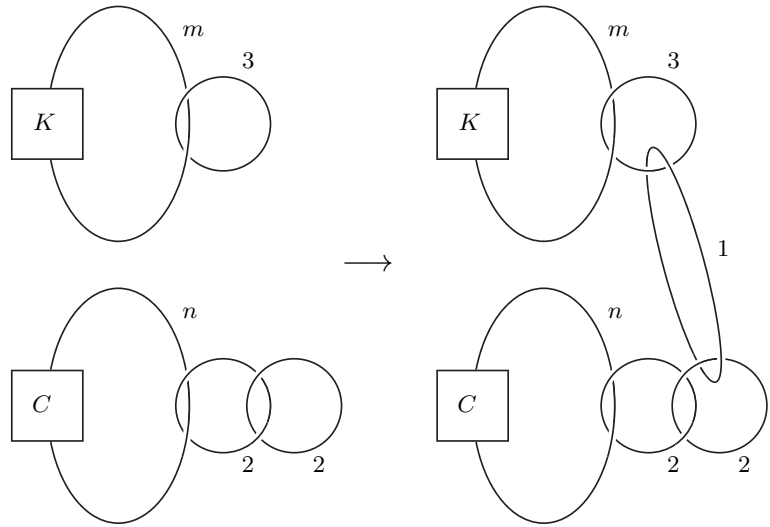}
\else \vskip 5cm \fi
\begin{narrow}{0.3in}{0.3in}
\caption{
\bf{A negative-definite cobordism from $S^3_{m-1/3}(K)\sqcup S^3_{n-2/3}(C)$ to $S^3_{m+n-1}(K\#C)$.}}
\label{fig:sum2}
\end{narrow}
\end{center}
\end{figure}

If $r=m-q/p$ and $s=n-(p-q)/p$ with $m,n $ positive integers and $\dfrac pq=[c_1, c_2, \dots, c_n]^+$ we may modify the diagram in Figure \ref{fig:sum1} accordingly.  
See Figure \ref{fig:sum2} for the case $q/p=1/3$. 
On the left hand side, 
we have $K$ with framing $m$ and a chain of unknots with framings 
$$c_1+1,\underbrace{2,\ldots,2}_{c_2-1},c_3+2,\underbrace{2,\ldots,2}_{c_4-1},c_5+2,\ldots$$
 and also $C$ with framing $n$ and a chain of unknots with framings 
$$\underbrace{2,\ldots,2}_{c_1-1},c_2+2,\underbrace{2,\ldots,2}_{c_3-1},c_4+2,\ldots.$$
By Lemmas \ref{lem:cf1} and \ref{lem:cf2} this is a four-manifold bounded by  $S^3_r(K)\sqcup S^3_s(C)$.
The cobordism $W$ is obtained by adding a single $+1$ framed unknot which links each of the two rightmost unknots described above once.
A sequence of $(+1)$-blowdowns converts this diagram to the right hand diagram in Figure \ref{fig:sum1} (with one of $m,n$ decreased by 1), from which it follows that $W$ is again a negative-definite cobordism
from $S^3_r(K)\sqcup S^3_s(C)$ to $S^3_{r+s}(K\#C)$.
\endproof

\begin{figure}[htbp] 
\begin{center} 
\ifpic
\psfrag{K}{\scriptsize$K$}
\psfrag{C}{\scriptsize$C$}
\psfrag{m}{\scriptsize$m$}
\psfrag{n}{\scriptsize$n$}
\psfrag{0}{\scriptsize$0$}
\psfrag{1}{\scriptsize$-1$}
\psfrag{2}{\scriptsize$-2$}
\psfrag{4}{\scriptsize$-4$}
\psfrag{arrow}{$\longrightarrow$}
\psfrag{boundary}{$\ds\stackrel{\partial}{\sim}$}
\includegraphics[width=12cm]{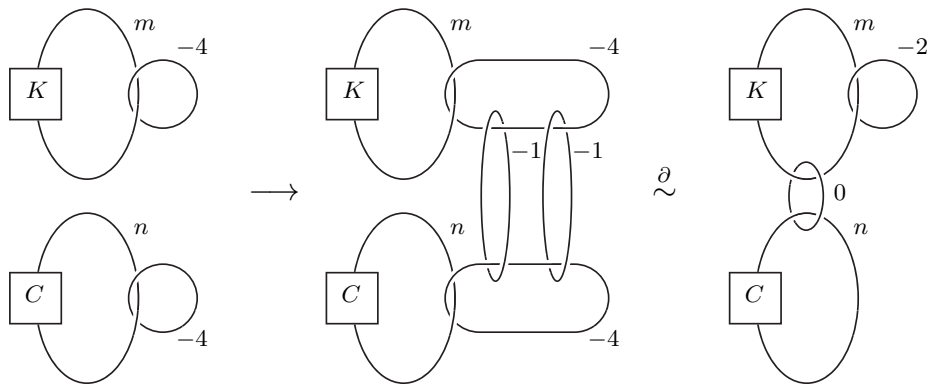}
\else \vskip 5cm \fi
\begin{narrow}{0.3in}{0.3in}
\caption{
\bf{A negative-definite cobordism from $S^3_{m+1/4}(K)\sqcup S^3_{n+1/4}(C)$ to $S^3_{m+n+1/2}(K\#C)$.}}
\label{fig:sum3}
\end{narrow}
\end{center}
\end{figure}

\begin{lemma}
\label{lem:sum2}
Let $K$ and $C$ be knots in $S^3$ and let $l,m,n$ be nonnegative integers with $l>1$.  Then there exists a negative-definite cobordism from $S^3_{m+1/2l}(K)\sqcup S^3_{n+1/2l}(C)$ to $S^3_{m+n+1/l}(K\#C)$.
\end{lemma}
\proof
The cobordism is illustrated in Figure \ref{fig:sum3} for the case $l=2$.  On the left hand side we have a Kirby diagram representing a four-manifold with boundary
$S^3_{m+1/2l}(K)\sqcup S^3_{n+1/2l}(C)$.  The cobordism $W$ is given by adding $l$ $(-1)$-framed 2-handles as shown.  This preserves $b_2^+$ and increases $b_2^-$ by $l$ so $W$ is negative-definite.  Blowing down the $(-1)$-framed 2-handles and sliding one of the resulting $l$-framed 2-handles over the other results in the last diagram shown in Figure \ref{fig:sum3}; as in Lemma \ref{lem:sum1} this is seen to represent the three-manifold $S^3_{m+n+1/l}(K\#C)$.
\endproof


\section{Basic properties}
\label{sec:basics}

In this section we establish some properties of $m(K)$, in particular its existence.

\proof[Proof of Theorem \ref{thm:properties}]
(a) Let $Y$ be the $4p$-surgery on $K$. If all the crossing changes in the unknotting of $K$ are realized by $(-1)$-blow-ups, the resulting surgery description of $Y$ consists of a $0$-framed unknot (corresponding to $K$) along with an unlink of $(-1)$-framed unknots. Hence $Y$ bounds a negative-definite $4$-manifold with one $1$-handle. 

(b) Given $r=p/q>m(K)$ there exists a rational number $s=p'/q'\in(m(K),r)$ with $p,p'$ coprime. 
By the definition of $m(K)$ and Lemma \ref{lem:rscob}, $S^3_s(K)$ bounds a negative-definite four-manifold $X_s$. Let $W$ be the cobordism from $S^3_s(K)$ to $S^3_r(K)$ given by Lemma \ref{lem:rscob}. Since $W$ is a two-handle cobordism, its first homology is a quotient of the first homology of either boundary component.
Since the orders of the first homology groups of the boundary components are coprime, they have no nontrivial common quotient and so $H_1(W)=0$. Taking the union $X_s \cup W$ yields a negative-definite manifold $X_r$ bounded by $S^3_r(K)$ with the property that the inclusion of the boundary in the manifold induces trivial homomorphism on $H_1$. Using Poincar\'e duality this implies that the restriction map on $H^2$ is onto. 

(c) If $K'$ is concordant to $K$, let $A \subset S^3 \times I$ be the annulus realizing the concordance. Then $r$-surgeries on $K$ and $K'$ extend over $A$; the resulting $4$-manifold is a homology cobordism from $S^3_r(K)$ to $S^3_r(K')$.

(d) Let $m=\ceil{m(K)}$ and $n=\ceil{m(C)}$. Suppose first that both $S^3_m(K)$  and $S^3_n(C)$ bound negative-definite four-manifolds.  Gluing these to the cobordism from Lemma \ref{lem:sum1} yields a negative-definite manifold bounded by
$S^3_{m+n}(K\#C)$ so that $\ceil{m(K\#C)}\le m+n$.

In  general  $S^3_r(K)$  and $S^3_s(C)$ bound negative-definite four-manifolds for any rational surgery coefficients $r>m$ and $s>n$.  In particular we may take
$r=m+1/2l$ and $s=n+1/2l$ for any positive integer $l$.  Combining with the negative-definite cobordism from Lemma \ref{lem:sum2} we see that $S^3_{m+n+1/l}(K\#C)$ bounds negative-definite.  Letting $l\to\infty$ we again see that $\ceil{m(K\#C)}\le m+n$.
\endproof


\section{Torus knots}
\label{sec:torus}
In this section we prove Theorem \ref{thm:Tpq} and Corollary \ref{cor:Tpq}.

Let $p>q>1$  be coprime integers and let $\frac pq=[c_1,c_2,\ldots,c_n]^+$, $c_i>0$ and $c_n\ge 2$.   Let
$$\mu(p,q)=
\left\{
\begin{array}{ll}
\ds pq-\frac{q}{p^*} & \quad\mbox{if $n$ is even,}\\
\ds pq-\frac{p}{q^*} & \quad\mbox{if $n$ is odd.}\\
\end{array}
\right.$$

\begin{proposition}\label{prop:Tpqobstr}
If $S^3_r(\Tpq)$ bounds a negative-definite four-manifold then $r\ge\mu(p,q)$.
\end{proposition}

\begin{proposition}\label{prop:Tpqemb}
The manifold $S^3_{\mu(p,q)}(\Tpq)$ embeds in a connected sum of $\cp$'s as a separating submanifold, and hence bounds a negative-definite four-manifold.
\end{proposition}

\begin{proposition}\label{prop:Tpqtors}
If $W$ is any negative-definite manifold that $S^3_{\mu(p,q)}(\Tpq)$ bounds then the restriction homomorphism $H^2(W;\zz)\to H^2(S^3_{\mu(p,q)}(\Tpq);\zz)$ is not onto; consequently, $H_1(W;\zz)$ contains nontrivial torsion.
\end{proposition}

\begin{figure}[htbp]
\begin{center}
\psfrag{e}{\scriptsize$e$}
\psfrag{f 1}{\scriptsize$\displaystyle\frac{\alpha_1}{\beta_1}$}
\psfrag{f 2}{\scriptsize$\displaystyle\frac{\alpha_2}{\beta_2}$}
\psfrag{f 3}{\scriptsize$\displaystyle\frac{\alpha_3}{\beta_3}$}
\includegraphics[width=10cm]{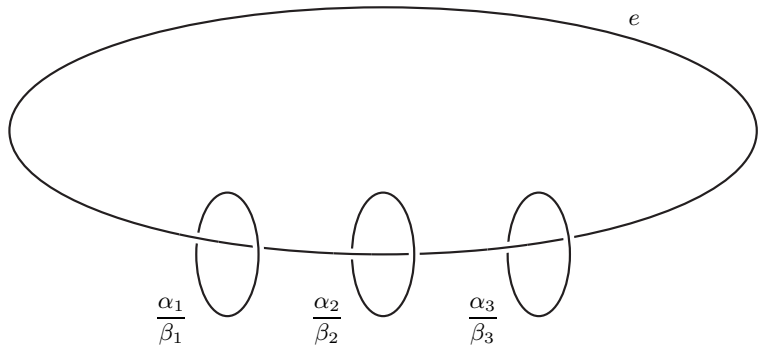}
\caption{$Y(e;\frac{\alpha_1}{\beta_1},\frac{\alpha_2}{\beta_2},\frac{\alpha_3}{\beta_3})$.}\label{fig:SFS}
\end{center}
\end{figure}

We use notation $Y(e;\frac{\alpha_1}{\beta_1},
\frac{\alpha_2}{\beta_2},\frac{\alpha_3}{\beta_3})$ 
to denote the 3-manifold that results by performing surgeries with the listed fractional coefficients on disjoint fibres of the degree $e$ $S^1$-bundle over $S^2$, as in Figure \ref{fig:SFS}. 
If the fractional coefficients are nonzero this is a Seifert fibred space whose exceptional fibres have orders 
$\alpha_i$.  We will also allow $\frac{\alpha_i}{\beta_i}$ to be zero or $\infty$.

\begin{lemma}[cf. Moser \cite{moser}]\label{lem:SFS}
For any rational number $r$,
$$S^3_r(\Tpq)=Y\left(2;\frac{p}{q^*},\frac{q}{p^*},\frac{pq-r}{pq-r-1}\right),$$
\end{lemma}

\proof
Start with the Seifert fibration of $S^3$ by $(p,q)$ torus knots which has two exceptional orbits, one of order $p$ and the other of order $q$. To determine the surgery coefficient (relative to the fibration) corresponding to the $r$-surgery along a regular fibre note that the linking number $\lk(K,K')=pq$, where $K$ is the surgery curve and $K'$ is a nearby regular fibre. For $r=a/b$ the surgery curve is given by $\gamma=a\mu+b\lambda=(a-bpq)\mu+bK'$, which yields $r-pq$ for the surgery coefficient.  It follows that
$$S^3_r(\Tpq)=Y\left(0;\frac{p}{\beta_1},\frac{q}{\beta_2},r-pq\right),$$
for some $\beta_1$, $\beta_2$.  Then the order of the first homology of the surgery is
$$\beta_1q(a-pqb)+\beta_2p(a-pqb)+pqb=\pm a \implies$$
$$pqb(q\beta_1+p\beta_2-1)=a(q\beta_1+p\beta_2\mp 1).$$
Since the last equation holds for all $r=a/b$ we conclude $q\beta_1+p\beta_2-1=0$.

We can fix coefficients $\beta_1$ and $\beta_2$ by requiring $0<\beta_1<p$ and $|\beta_2|<q$, so that $\beta_1=q^*$
and $\beta_2+q=p^*$.  The result now follows by applying Rolfsen twists (see for example \cite[\S 5.3]{gs}) 
to the unknots in 
the Dehn surgery diagram for $Y\left(0;p/\beta_1,q/\beta_2,r-pq\right)$ with the framings $q/\beta_2$ and $r-pq$.
\endproof

\proof[Proof of Proposition \ref{prop:Tpqobstr}.]
We may assume that $r<pq-1$.   Using Lemmas \ref{lem:SFS}, \ref{lem:cf0}, \ref{lem:cf1} and \ref{lem:cf3}   we find that
$S^3_r(\Tpq)$ is the boundary of a positive-definite plumbing $P$ of disk bundles over spheres 
corresponding to a tree with $3$ legs where the weight of the central vertex is $2$.  The weights on the three legs (listed from the central vertex) are as follows:
\begin{itemize}
\item the weights on the first leg are the coefficients in the negative continued fraction for $\frac{pq-r}{pq-r-1}$;
\item the weights on the second leg are
$$ c_n+1,\underbrace{2,\ldots,2}_{c_{n-1}-1},c_{n-2}+2, \ldots, \underbrace{2,\ldots,2}_{c_4-1},c_3+2,\underbrace{2,\ldots,2}_{c_2-1},c_1+1$$
if $n$ is odd and
$$ c_n+1,\underbrace{2,\ldots,2}_{c_{n-1}-1},c_{n-2}+2,\ldots, \underbrace{2,\ldots,2}_{c_5-1},c_4+2,\underbrace{2,\ldots,2}_{c_3-1},c_2+1$$
if $n>2$ is even (and $c_2$ if $n=2$).
\item the weights on the third leg are
$$\underbrace{2,\ldots,2}_{c_{n}-1},c_{n-1}+2, \ldots, \underbrace{2,\ldots,2}_{c_5-1},c_4+2,\underbrace{2,\ldots,2}_{c_3-1},c_2+1$$
if $n$ is odd and
$$ \underbrace{2,\ldots,2}_{c_{n}-1},c_{n-1}+2,\ldots, \underbrace{2,\ldots,2}_{c_4-1},c_3+2,\underbrace{2,\ldots,2}_{c_2-1},c_1+1$$
if $n$ is even.
\end{itemize}

(Note that if $n$ is odd, then the second leg arises from the continued fraction expansion of $\ds\frac p{q^*}$ and the third leg corresponds to $\ds\frac q{p^*}$; for $n$  even $\ds\frac q{p^*}$ gives the second leg and $\ds\frac p{q^*}$ gives the third.)

If $S^3_r(\Tpq)$ bounds a negative-definite 4-manifold $X$, then $P\cup (-X)$ is a closed positive-definite manifold so by Donaldson's Theorem the intersection form of $P$ embeds in some $\zz^k$ (with the standard form).  We seek the minimal (or a priori, infimal) $r$ for which such an embedding is possible.  We first note that minimising $r$ is equivalent to minimising 
$$\frac{pq-r}{pq-r-1}=[a_1,\dots,a_m]^-=[a_1,\dots,a_m,\infty]^-,$$
which in turn is equivalent to 
 finding the smallest sequence $a_1,a_2,\ldots,a_m,\infty$ in lexicographical ordering, with integer coefficients $a_i\ge2$.

Let $E$ denote the central vertex of $P$, and let $U_1, U_2,\dots$ denote the vertices on the first leg.  Similarly label the vertices on the second and third legs with $V_i$ and $W_j$ respectively.
Denote basis vectors of $\zz^k$ by $e_i$ and $f_j$.  Suppose for some $r$ we have an embedding of the intersection form of $P$ in $\zz^k$.  Without loss of generality, $E$ maps to $e_1+f_1$ and $V_1$ maps to $-e_1+x$ for some $x\in\zz^k$.

We claim there is an embedding with $a_1=2$ so that $U_1$ maps to one of $-e_1+e_2$ or $-f_1+f_2$.  Suppose first that $U_1$ maps to $-e_1+e_2$, so that $V_1$ maps to $-e_1-e_2+x'$.
If we can take $a_2=2$ as well we must map $U_2$ to $-e_2+e_3$ and $V_1$ to $-e_1-e_2-e_3+x''$.  Continuing in this way we find that $U_i$ maps to $-e_i+e_{i+1}$ for $i=1,\dots,c_n-1$, $V_1$ maps to $e_1+e_2+\dots+e_{c_n+1}$ and $V_2$ maps to $e_{c_n+1}-e_{c_n+2}$, and so on.  The requirement that $r$ be minimal combined with the assumption that $U_1$ maps to $-e_1+e_2$ completely determines the weights on the first leg and the embedding in $\zz^k$ (up to automorphism of $\zz^k$) of the first two legs.  The weights on the first leg (under this assumption) are
$$\underbrace{2,\ldots,2}_{c_{n}-1},c_{n-1}+2, \ldots,\underbrace{2,\ldots,2}_{c_3-1},c_2+2,\underbrace{2,\ldots,2}_{c_1-1},$$
if $n$ is odd and
$$ \underbrace{2,\ldots,2}_{c_{n}-1},c_{n-1}+2,\ldots, \underbrace{2,\ldots,2}_{c_4-1},c_3+2,\underbrace{2,\ldots,2}_{c_2-1},$$
if $n$ is even.

The reasoning from the previous paragraph may be applied to the third leg instead of the first, showing that $W_1$ maps to $-f_1+f_2$ (the weights on the third leg represent a smaller continued fraction than the minimal value found in the previous paragraph).  This enables us to eliminate the possibility that $U_1$ maps to $-f_1+y$ for any $y\in\zz^k$ since orthogonality with the vertices on the third leg would then imply $a_1>2$.

It remains to extract the value of the minimal $r$.  We have
$$\frac {pq-r}{pq-r-1}=[\underbrace{2,\ldots,2}_{c_{n}-1},c_{n-1}+2, \ldots,\underbrace{2,\ldots,2}_{c_{m+2}-1},c_{m+1}+2,\underbrace{2,\ldots,2}_{c_m-1}]^-,$$
where $m=1$ if $n$ is odd, else $m=2$.  Applying (\ref{eqn:cf}) and Lemmas \ref{lem:cf0}, \ref{lem:cf1} and \ref{lem:cf3} we find
\begin{eqnarray*}
pq-r&=&[c_n,c_{n-1},\dots,c_m]^+\\
&=&
\left\{
\begin{array}{ll}
\ds \frac{q}{p^*} & \quad\mbox{if $n$ is even,}\\
\ds \frac{p}{q^*} & \quad\mbox{if $n$ is odd.}\\
\end{array}
\right.
\end{eqnarray*}

\endproof

\proof[Proof of Proposition \ref{prop:Tpqemb}.]
Assume for convenience that $n$ is odd (the even case is proved in exactly the same way).
From the proof of Proposition \ref{prop:Tpqobstr} we see that $S^3_{\mu(p,q)}$ is the boundary of the positive-definite 4-manifold $P$ presented by a Kirby diagram corresponding to a three-legged tree, with framing $2$ on the unknot corresponding to the central vertex, and framings on the 3 legs given by
\begin{itemize}
\item 
first leg: $\underbrace{2,\ldots,2}_{c_{n}-1},c_{n-1}+2, \ldots,\underbrace{2,\ldots,2}_{c_3-1},c_2+2,\underbrace{2,\ldots,2}_{c_1-1},$
\item 
second leg: $ c_n+1,\underbrace{2,\ldots,2}_{c_{n-1}-1},c_{n-2}+2, \ldots, \underbrace{2,\ldots,2}_{c_4-1},c_3+2,\underbrace{2,\ldots,2}_{c_2-1},c_1+1$,
\item 
third leg: $\underbrace{2,\ldots,2}_{c_{n}-1},c_{n-1}+2, \ldots, \underbrace{2,\ldots,2}_{c_5-1},c_4+2,\underbrace{2,\ldots,2}_{c_3-1},c_2+1$.
\end{itemize}

We add further handles to this diagram to get a diagram for a manifold $\widetilde{X}$.  For each vertex in the third leg, if the weight of the vertex is $w$ and the valency is $v$ (one if the rightmost vertex, or else 2), add $w-v$ parallel $(+1)$-framed meridians to the corresponding component of the Kirby diagram. We claim the resulting $\widetilde{X}$ can be obtained from $S^2\times B^2$ by a sequence of $(+1)$-blow-ups.

Begin by blowing down $(+1)$-framed unknots on the third leg; this completely eliminates the third leg and replaces the weight on the central vertex by $1$.  We now have a linear plumbing with weights
$$c_1+1,\underbrace{2,\ldots,2}_{c_2-1},c_3+2,\dots,\underbrace{2,\ldots,2}_{c_{n-1}-1},c_n+1,1,
\underbrace{2,\ldots,2}_{c_{n}-1},c_{n-1}+2, \ldots,\underbrace{2,\ldots,2}_{c_3-1},c_2+2,\underbrace{2,\ldots,2}_{c_1-1}.$$
A simple induction argument shows that successive blow-downs reduce this to a single zero-framed unknot. Hence we can add a 3-handle and 4-handle to $\widetilde{X}$ to get a connected sum of $\cp$'s.

Let $X$ be the closure of the complement of $P$ in the connected sum of $\cp$'s; then $-X$ is a negative-definite manifold bounded by 
$S^3_{\mu(p,q)}(T_{p,q})$.
\endproof

A sublattice $L\subset \zz^k$ is called {\em primitive} if the quotient $\zz^k/L$ is torsion-free. If $L$ is not primitive in $\zz^k$, then the restriction homomorphism on the dual lattices $\Hom(\zz^k,\zz)\to L'=\Hom(L,\zz)$ is not surjective. This observation yields the following result which is the main ingredient in the proof of Proposition \ref{prop:Tpqtors}.

\begin{lemma}\label{lem:torsion}
Let $Y$ be a rational homology $3$-sphere that bounds a positive-definite $4$-manifold $X$ with $H_1(X;\zz)=0$.  
Suppose that the intersection lattice $L=(H_2(X;\zz),Q_X)$ of $X$ does not admit a primitive embedding in any $\zz^k$.
Then for any negative-definite $4$-manifold $W$ that $Y$ bounds the restriction homomorphism $H^2(W;\zz)\to H^2(Y;\zz)$ is not onto; consequently, $H_1(W;\zz)$ contains nontrivial torsion.
\end{lemma}

\proof 
It follows from the exact cohomology sequence of the pair $(X,Y)$ (using assumption $H_1(X;\zz)=0$) 
that $H^2(Y;\zz)\cong L'/L$, where $L'$ denotes the dual lattice of $L$. Let $W$ be a negative-definite manifold with boundary $Y$; we may assume that $b_1(W)=0$. Let $Z=X \cup_Y (-W)$. It follows from the Mayer-Vietoris homology exact sequence for $Z$ that $L$ embeds in $H_2(Z)/\mathrm{Tors}$. Note that the intersection pairing of $Z$ is isomorphic to some $\zz^k$ by Donaldson's Theorem. Since by assumption $L$ is not primitive in $\zz^k$, the restriction $H^2(Z) \to L'$ is not onto. 
Consider now the Mayer-Vietoris cohomology exact sequence for $Z$:
$$\cdots\to H^2(Z) \to H^2(X)\oplus H^2(W) \to H^2(Y) \to \cdots .$$
Choose an element $x\in L'= H^2(X)$ that is not in the image of the restriction homomorphism from $H^2(Z)$. Then the image of $x$ in $H^2(Y)$ is not in the image of the restriction $H^2(W)\to H^2(Y)$; if it were the image of some $y\in H^2(W)$ then $x\oplus y$ would be in the image of $H^2(Z)$.
\endproof

\proof[Proof of Proposition \ref{prop:Tpqtors}]
We again assume for clarity of exposition that $n$ is odd.  The proof of Proposition \ref{prop:Tpqobstr}
shows that $S^3_{\mu(p,q)}(T_{p,q})$ bounds a plumbing $P$ whose intersection form (call it $L$) embeds in $\zz^k$.  We claim that for any such embedding $L$ is not primitive in $\zz^k$.  Proposition \ref{prop:Tpqtors} is a consequence of this claim and Lemma \ref{lem:torsion}.

From the proof of Proposition \ref{prop:Tpqobstr} we see that there is a unique embedding in $\zz^k$ (up to automorphisms) of the first two legs of the tree defining $P$.  A simple recursive description of all such embeddings is as follows: starting with the sequence $-e_1-e_2,e_2,e_1-e_2$ one applies a finite sequence of the following modifications (blow-ups):
$$\dots,v,e_i,w,\ldots\quad \leadsto\quad \dots,v-e_{i+1},e_{i+1},e_i-e_{i+1},w,\dots$$
or
$$\dots,v,e_i,w,\ldots\quad \leadsto\quad \dots,v,e_i-e_{i+1},e_{i+1},w-e_{i+1},\dots$$
and then replaces the final $e_s$ with $e_s+f_1$.  (This becomes the central vertex of the tree, and the two chains on either side give the first two legs of the tree; comparing to the proof of Proposition \ref{prop:Tpqobstr} one should reverse the order of indices of the basis vectors $e_1,\dots,e_s$.)  
Inductively we see that the image of the sublattice $L_0$ of $L$ corresponding to the first two legs rationally spans $\zz^{s}$ but is not equal to it as its determinant is $p^2>1$. We claim that the image of $L$ intersected with this $\zz^s$ is equal to the image of $L_0$ and therefore $L$ is not primitive in $\zz^k$. Indeed, the embedding of the third leg cannot use any of the basis vectors $e_1,\dots,e_s$, hence the only way the intersection could be larger is if some multiple of $e_s$ were in the image of the central vertex and the third leg. Using that $e_s$ is rationally in the image of $L_0$, so
$$e_s=\sum_i c_iv_i+\sum_j d_jw_j,$$
where $v_i$ ($w_j$) are images of vectors in the first (second) leg and $c_i,d_j\in \qq$, this would imply that the two sums above (representing orthogonal vectors) vanish, providing a contradiction.
\endproof

Theorem \ref{thm:Tpq} follows from Propositions \ref{prop:Tpqobstr}, \ref{prop:Tpqemb},  \ref{prop:Tpqtors}, and Theorem \ref{thm:properties}(a).

\proof[Proof of Corollary \ref{cor:Tpq}]
The proof is based on that of \cite[Theorem 4.2]{ls}; we briefly recall the argument here.  Let $Y_r$ denote the surgery manifold $S^3_r(T_{p,q})$.
By \cite[Proposition 4.1]{ls}, $Y_r$ is an $L$-space whenever
$r\ge2g_s(T_{p,q})-1=pq-p-q$.  By \cite[Theorem 1.4]{gbounds}, any symplectic filling of an $L$-space is negative-definite.  By Theorem \ref{thm:Tpq}, $Y_r$ does not bound 
a negative-definite four-manifold if $r\in[pq-p-q,m(T_{p,q}))$ and thus does not admit a fillable contact structure.
%
%
%
%
\endproof

We note that the result in Corollary \ref{cor:Tpq}  is optimal: for any $r\notin [pq-p-q,m(T_{p,q}))$, the three-manifold obtained by $r$ surgery on $T_{p,q}$ \emph{does} admit a fillable contact structure.  This may be deduced using the classification by Lecuona-Lisca of Seifert fibred spaces which admit fillable contact structures \cite[Theorem 1.3]{ll}.



\begin{thebibliography}{99}
\bibitem{d} S.~K.~Donaldson.  \textsl{An application of gauge
theory to four-dimensional topology}, J.~Diff.~Geom. {\bf 18}
(1983), 279--315.
\bibitem{gs} R.~E.~Gompf \& A.~I.~Stipsicz. \textsl{{$4$}-manifolds and {K}irby calculus},
      Graduate \allowbreak Studies in Math. {\bf 20}, Amer.~Math.~Soc., 1999.
\bibitem{greene} J.~E.~Greene. \textsl{L-space surgeries, genus bounds, and the cabling conjecture}, arXiv:1009.1130, 2010.
\bibitem{hnk} F.~Hirzebruch, W.~D.~Neumann \& S.~S.~Koh. \textsl{Differentiable
    manifolds and quadratic forms},
    Lecture Notes in Pure and Applied Math. {\bf 4}, Marcel Dekker, 1971.
\bibitem{ll} A.~Lecuona \& P.~Lisca. \textsl{Stein fillable Seifert fibered 3-manifolds},  arXiv:1007.3010.
\bibitem{ls} P.~Lisca \& A.~I.~Stipsicz. \textsl{\ozsvath-\szabo\ invariants
and tight contact three-manifolds,~I}, Geom. Topol. {\bf 8} (2004), 925--945.
\bibitem{moser} L.~Moser. \textsl{Elementary surgery along a torus knot},
    Pacific Journal of Math. {\bf 38} (1971), 737--745.
\bibitem{neumann} W.~Neumann. \textsl{A calculus for plumbing applied to the topology of complex surface singularities and degenerating complex curves},
Trans. Amer. Math. Soc. {\bf 268} (1981), 299Ð344
\bibitem{short}B.~Owens \& S.~Strle.  \textsl{A characterisation of the $Z^n \oplus Z(\delta)$ lattice and definite nonunimodular intersection forms}, arXiv:0802.1495.
\bibitem{gbounds} P.~\ozsvath~ \& Z.~\szabo. \textsl{Holomorphic disks and genus bounds}, Geom. Topol. {\bf 8} (2004), 311--334.
\bibitem{ppp} P.~Popescu-Pampu. \textsl{The geometry of continued fractions and the topology of surface singularities}, Singularities in geometry and topology 2004, 119--195, Adv. Stud. Pure Math. {\bf 46}, Math. Soc. Japan, 2007.
\end{thebibliography}
\end{document}